\documentclass[a4paper]{amsart}
\usepackage[english]{babel}
\usepackage{amssymb}
\usepackage{graphicx}
\usepackage{booktabs}
\usepackage{longtable}
\usepackage[latin1]{inputenc}

\newcommand{\Z}{\mathbb{Z}}

\newcommand{\Q}{\mathbb{Q}}

\newtheorem{theo}{Theorem}

\theoremstyle{definition}

\theoremstyle{remark}

\begin{document}
\title[Elliptic curves and torsion]{Elliptic curves with torsion groups $\Z/ 8 \Z $ and $\Z/ 2 \Z \times \Z/ 6 \Z$}
\author[A. Dujella]{Andrej Dujella}
\address{Department of Mathematics\\Faculty of Science\\University of Zagreb\\Bijeni{\v c}ka cesta 30, 10000 Zagreb, Croatia}
\email[A. Dujella]{duje@math.hr}
\author[M. Kazalicki]{Matija  Kazalicki}
\address{Department of Mathematics\\Faculty of Science\\University of Zagreb\\Bijeni{\v c}ka cesta 30, 10000 Zagreb, Croatia}
\email[M. Kazalicki]{matija.kazalicki@math.hr}
\author[J. C. Peral]{Juan Carlos Peral}
\address{Departamento de Matem\'aticas\\Universidad del Pa\'is Vasco\\Aptdo. 644, 48080 Bilbao, Spain}
\email[J. C. Peral]{juancarlos.peral@ehu.es}
\thanks{A. D. and M. K. were supported by the Croatian Science Foundation under the project no. IP-2018-01-1313 and
the QuantiXLie Center of Excellence, a project co-financed by the
Croatian Government and European Union through the European Regional Development Fund - the Competitiveness and Cohesion Operational Programme (Grant
KK.01.1.1.01.0004). J. C. P. was supported by the UPV/EHU grant EHU 10/05.
J. C. P dedicate this paper to Ireneo Peral in memoriam.}
\begin{abstract}
In this paper, we present details of seven elliptic curves over $\mathbb{Q}(u)$ with rank $2$ and torsion group
$\mathbb{Z}/8\mathbb{Z}$ and five curves over $\mathbb{Q}(u)$ with rank $2$ and torsion group
$\mathbb{Z}/2\mathbb{Z}\times \mathbb{Z}/6\mathbb{Z}$. We also exhibit some particular examples of curves
with high rank over $\mathbb{Q}$ by specialization of the parameter. We present several sets of infinitely many
elliptic curves in both torsion groups and rank at least $3$ parametrized by elliptic curves having positive rank.
In some of these sets we have performed calculations about the distribution of the root number.
This has relation with recent heuristics concerning the rank bound for elliptic curves by Park, Poonen, Voight and Wood.
 \end{abstract}

\maketitle

\section{Introduction }
We are interested in curves having torsion group $\Z/ 8 \Z $  or  $\Z/ 2 \Z \times \Z/ 6 \Z$. We describe a general model for curves with these  torsions,  which jointly with the curves having torsion $\Z/ 7 \Z $ are the three cases in which the general model is a $K3$ surface.  Then we try to impose the existence of new points.

\subsection{New points from quadratic conditions}

We explain how we get new points in our models.

1) All the  curves in the paper  can be written as $y^2= x^3+  A x^2+ B x$. This is due to the fact that all of them have at least one torsion point of order two.

2) In these cases it is well known that the $x$-coordinates of rational points   should be either divisors of $B$ or rational  squares times divisors of $B$.

3) We know some of the divisors of $B$, namely the polynomial factors of $B$.

4) We construct the list of all such  known divisors of $B$ and we look for those that become a new point after parametrizing a conic.

If $d$ is a divisor of $B$, we try to force $d$ to be the $x$-coordinate of a new point, so we have to impose that
$d^3+ A d^2+
B d= d^2( d + A + B/d)$ is a square. So for every known divisor of $B$ we consider the equation
\[
 d + A + \frac{B}{d}=\hbox{Square}.
\]
In some cases this equation is equivalent to solving a conic and in this way we get a new point in the curve in case of solvability.

We try the same with a divisor of $B$, say $d$,  times a rational square, say $\frac{U^2}{V^2}$,  and
we  get the following equations
\[
d^3\frac{U^6}{V^6}+ A d^2\frac{U^4}{V^4}+B d\frac{U^2}{V^2}=
\frac{d^2 U^2}{V^6}\left(d U^4+ A U^2 V^2+V^4 \frac{B}{d}\right)=\hbox{Square}.
\]
These equations describe the homogeneous spaces corresponding to the pair $(U,V)$.
In some cases this equation is equivalent to a quadratic condition and in case of solvability we get a new point
on the curve, which in case of independence raise the rank by one.

\subsection{Torsion group $\Z/ 8 \Z $}
For torsion  $\Z/ 8 \Z $,
  we list the main  data for eighteen curves with rank $1$.
Curves having    rank $1$ over $\Q(u)$ were found
by Kulesz and Lecacheux and Woo, see  \cite{Ku}, \cite{Le1}, \cite{Le2} and \cite{W}.

Kulesz used  a quadratic section while
Lecacheux used fibrations  of the corresponding surfaces, as the ones given explicitly  in Beauville \cite{Be}, in  Bertin and Lecacheux  \cite{BL} or Livn\'e and Yui  \cite{LY}, in order to find elliptic curves  with positive rank over $\Q(u)$ for several torsion groups.

Woo in his Doctoral Dissertation \cite{W}
studied  in a systematic way quadratic sections of $K3$ surfaces and in this way he found several quadratic sections  leading to rank $1$ curves over $\Q(u)$ for torsion groups $\Z/ 8 \Z $  and $\Z/ 2 \Z \times \Z/ 6 \Z$ in which we are interested in this paper. These results can be found in \cite[Chapter~3]{W}.

In 2013, Dujella and Peral  (see \cite{DP1} and  \cite{DP2}) described ten of these curves, including those previously discovered by Kulesz and Lecacheux and Woo.
The main details of these ten curves are given first.
In 2014, MacLeod (see \cite{ML}) listed twelve examples. Six of  these twelve  are new, and we list them in second place. Finally, we discovered two new examples wich are given in third place.
We list all such curves in order to complete the presentation.
That they are non-isomorphic is deduced by comparing  their $j$-invariants.

 In some cases, we are able to repeat  the procedure of quadratic sections,  and so we reach  curves with rank $2$, as mentioned before. Details will be the subject of another section.
As already said, the first two examples for such curves with rank $2$ over $\Q(u)$ were found by  Dujella-Peral.  Details were  presented  first on ArXiv \cite{DP1} and later in \cite{DP2}.
In 2014, MacLeod found other two curves with rank $2$,  see \cite{ML}.
In this paper, we construct three new curves with rank $2$ over $\Q(u)$ and several examples of new sets of infinitely many elliptic curves with such torsion group and rank at least $3$ parametrized by elliptic curves having positive rank.
Again an argument with their $j$-invariant shows that the seven curves are different.

We also exhibit particular examples of curves with high rank by specialization of the parameter.
In particular,  for the torsion group  $\Z/ 8 \Z $, the rank record for elliptic curves over $\Q$ at the moment is $6$. The first example  is due to Elkies (2006), the second was found by Dujella-MacLeod-Peral (2013)
and the third and fourth by Voznyy (2021).
See \cite{D} for details of these curves were also the details of  more than $50$  curves with rank $5$ are given.

\subsection{Torsion group $\Z/2\Z \times \Z/6\Z$ }
Elliptic curves with torsion group $\Z/2\Z \times \Z/6\Z$ and rank at least $1$ over $\Q(u)$ have been constructed by several authors,  see  \cite{Ca}, \cite{Ku}, \cite{Le1}, \cite{Du1} and \cite{Ra}, and as mentioned before by Woo in \cite{W}. For example the transformation  leading the Woo model (written in terms of $t$) to our model is $t\to -\frac{2 (v-1)}{v+1}$.

Here we describe five  elliptic curves with  torsion group
$\Z/2\Z \times \Z/6\Z$ and rank $2$ over $\Q(u)$.

The first four examples with rank $2$ over $\Q(u)$ and torsion $\Z/2\Z \times \Z/6\Z$  were found by  Dujella-Peral in  2013.
Three of them were found using Diophantine triples, and the fourth by starting with a model due to Hadano. The last example is new.
 We have checked that  the five  curves  are non-isomorphic to each other  by comparing their $j$-invariants.
 It should be observed  that the seven curves  presented by MacLeod in 2014 are all isomorphic to one of the first three curves.

Let us mention at this point that in \cite{DP3} it was shown that the elliptic curves over $\Q (t)$ induced by  Diophantine triples can have as torsion group any of the non cyclic groups in Mazur's theorem,
and several results with highest known ranks for such curves were obtained.
In \cite{DJS}, it was proved that there are infinitely many Diophantine triples over quadratic fields which induce elliptic curves with torsion groups $\Z/2\Z \times \Z/10\Z$,  $\Z/2\Z \times \Z/12\Z$ and $\Z/4\Z \times \Z/4\Z$.

 We also exhibit, in the present paper, some particular examples of curves with high rank over $\Q$ by specialization of the parameter.
In the case of torsion group  $\Z/2\Z \times \Z/6\Z$, the rank record for elliptic curves over $\Q$ at the moment is $6$. The first example was found by Elkies  in 2006, the second was found by Dujella-Peral-Tadi\'c in 2015
and the third and fourth by Dujella-Peral in 2020.  See \cite{D} for details of these curves where also the details of  more than $60$  curves with rank $5$ are given.

We include  several  sets of infinitely many elliptic curves in both  torsion groups having  rank at least $3$ parametrized by elliptic curves having positive rank.
 In some of these sets,  we have performed
 calculations concerning the distribution of the root number. This has relation with recent  heuristics about the rank bound for elliptic curves by Park, Poonen, Voight and Wood \cite{PPVW}.
 In their paper, the authors describe a heuristic argument according to which the rank for elliptic curves over $\Q$ would be bounded.
 For example, they assert that only a finite number of curves will have rank greater than $21$.
 And for the torsion groups considered here, only a finite number of curves will have rank greater than $3$.

 This is in contrast with an old conjecture predicting the existence of elliptic curves of arbitrarily high rank over $\Q$ even for each torsion group in the Mazur theorem.
 Now observe that if both the parity conjecture and the predictions of the heuristic were true,
 then the distribution of root number over these infinite  families with rank $3$ would not be evenly distributed.
 We  dedicate a section for the  calculation  related to these questions.

\section{General model for  torsion groups $\Z/8\Z $ and $\Z/2\Z \times \Z/6\Z$. }
The Tate normal form for an elliptic curve is given by
\[
E(b,c): \quad y^2+(1-c) x y - b y = x^3- b x^2
\]
(see \cite{Kn}). It is nonsingular if and only if $b\neq 0$.

  Using the addition law for $P=(0,0)\in E(b,c)$  and taking $d=b/c$  (assuming that $c\neq 0$),
  we have, in terms of $c$ and $d$:
\begin{align*}
P &=(0,0),&
-P &=(0, c\, d)\\
 2 P&=(c \,d,0),&
- 2 P&=(0,b)\\
 3 P&=(c,c\,(d-1)),&
- 3 P&=(c,c^2)\\
4 P&=\left( d (d -1), d^2 ( c- d+1)\right),&
- 4 P&=\left(d (d-1), d(d-1)^2\right).
\end{align*}
From here we get that   $P$ is a torsion point  of order $6$ for $d=1+c$ or  $b=c+c^2$, $P$ is a torsion point of order $7$ for
$c=d(d-1)$ and $b=d^2(d-1)$ and
$P$ is a torsion point  of order $8$ for
\begin{align*}
b&= (2 d-1)(d-1),\\
c&=\frac{(2 d-1)(d-1)}{d}
\end{align*}
with $d$   rational, see \cite{Kn}.

For these values  we can write the general curve  with torsion $\Z/ 8\Z$ in the form
$y^2= x^3+ A_{8}(v) x^2 + B_{8}(v) x$,  where
\begin{align*}
A_{8}(v)&=1 - 8 v + 16 v^2 - 16 v^3 + 8 v^4, \\
B_{8}(v)&=16 (v-1)^4 v^4
\end{align*}
(we replaced here $d$ by $v$ to adjust notation with $\Z/2\Z \times \Z/6\Z$ case
since we will similarly treat these two cases).

The $j$-invariant for this model is
\[
j_8(v)=\frac{\left(16 v^8-64 v^7+224 v^6-448 v^5+480 v^4-288 v^3+96 v^2-16 v+1\right)^3}{(v-1)^8 v^8 (2 v-1)^4
   \left(8 v^2-8 v+1\right)}.
\]

This invariant remains  the same with the substitution $v\to 1-v$.


Finally  for  torsion $\Z/ 6\Z$ the general curve is  $y^2= x^3+ A_{6}(d) x^2 + B_{6}(d) x$,  where
\begin{align*}
A_{6}(d)&=1 + 6 d- 3 d^2, \\
B_{6}(d)&=-16 d^3.
\end{align*}
Now the discriminant of $X^3+A_6 (d)X^2+ B_6 (d) X$ is $256 d^6 (d+1)^3 (9 d+1)$
so we get complete factorization, and hence torsion $\Z/2\Z \times \Z/6\Z$, if this is a square.
By taking $(d+1)(9d+1)=(3d+v)^2$, we get that the general curve with this torsion group can be written as:
$y^2= x^3+ A_{26}(v) x^2 + B_{26}(v) x$  where
\begin{align*}
A_{26}(v)&=37 - 84 v + 102 v^2 - 36 v^3 - 3 v^4, \\
B_{26}(v)&=32 (v-1)^3 (v+1)^3 (3 v-5).
\end{align*}

The $j$-invariant for this model is

\[
j_{2 \times 6}(v)=\frac{\left(3 v^2-6 v+7\right)^3 \left(3 v^6-18 v^5+345 v^4-1260 v^3+1605 v^2-738 v+127\right)^3}{4 (v-3)^6
   (v-1)^6 (v+1)^6 (3 v-5)^2 (3 v-1)^2}.
   \]

In this case the transformations leading to unchanged $j$-invariant are the following:
\[
\{v\to 2-v\},\left\{v\to \frac{v-7}{3 v-5}\right\},\left\{v\to \frac{v+5}{3 v-1}\right\},\left\{v\to
   \frac{5 v-7}{3 v-1}\right\},\left\{v\to \frac{5 v-3}{3 v-5}\right\}.
\]

%
%
%


As said before  in both  cases $\Z/ 8\Z$  and  $\Z/2\Z \times \Z/6\Z$,  this form of the curve, i.e., $y^2=x^3+A x^2+ B x$,  is a convenient way  to search for candidates for new rational points. In fact,  their $x$-coordinates  should be either divisors of $B$ or rational  squares times divisors of $B$.
For example, if we define
\[
el_8(x)=x^3+A_8(v) x^2+ B_8(v) x,
\]
 we look for values of $x$ such that the corresponding expression can be transformed into a square
 by solving a quadratic equation in $v$.
Observe that $el_8(4v^4)=16 v^8 (2 v-1)^2 (4 v^2-4 v+5)$,
so we get a new point by solving \[4 v^2-4 v+5=\hbox {Square.}\]
This is the first example in our list of rank $1$ curves with such torsion.

\section{Rank $1$   curves for  torsion group $\Z/8\Z $.}

We list eighteen points together with the condition and the specialization of  $v$, which leads to a new point on the curve and so, eventually,  to rank  $1$ curves over $\Q(w)$.
Observe that the parametrization it is not unique.

We split the list into three parts, having  ten, six and two entries, respectively.
First we present the ten conditions that were in the papers of Dujella and Peral already quoted.
As said, all of them were included in Dujella-Peral \cite{DP1} and \cite{DP2}. Some typos there are now  corrected.
The curve corresponding to place $4$ were already in Lecacheux, and  the curve in place $5$ is both in Lecacheux and  Kulesz.

In the second part  of the list, we present  the six new curves discovered by MacLeod in \cite{ML}.
MacLeod reported twelve curves of rank $1$  in his  paper, of which six are equivalent to some in the preceding ten and the other six  were new.  These curves form  the second part of the list. Finally, the last two entries correspond to new curves and are the third part of the list.

\subsection{First part of the list.}
\mbox{}

{\footnotesize
\begin{align*}
x_1 &=4 v^4,&c_1= 5 - 4 v + 4 v^2,&   &v_1=\frac{5-w^2}{4 (w+1)},\\
\\
x_2 &=-(-1 + v) v,&c_2=1+v-v^2,&   &v_2=\frac{(w-2) w}{w^2+1},\\
\\
x_3 &=-4 v^3 (-2+3 v),& c_3= - (2 + v) (-2 + 3 v),& &v_3=-\frac{2 (w-1) (w+1)}{w^2+3},\\
\\
x_4 &=16 (-1 + v)^2 v^2 (1 - 2 v + 2 v^2),&c_4=(1 - 2 v + 2 v^2),&   &v_4=\frac{(w-2) w}{w^2-2},\\
\\
x_5 &=-2 v^2 (-1+2 v^2),& c_5=- (-1 + 2 v^2) ,&  &v_5=-\frac{2 w}{w^2+2},\\
\\
x_6 &=-(-1 + v)^2 (1 - 6 v + 4 v^2),&c_6=-1 + 6 v - 4 v^2,&   &v_6=\frac{w^2-2 w+2}{w^2+4},\\
\\
x_7 &=-\frac{4 (v-1)^4 (6 v-1)}{2 v-3},&c_7=-(-3 + 2 v) (-1 + 6 v),&     &v_7=\frac{3 w^2+1}{2 \left(w^2+3\right)},\\
\\
x_8 &=-\frac{4 (v-1)^4 (4 v-1)}{4 v-3},&c_8=- (-3 + 4 v) (-1 + 4 v), &    &v_8=\frac{w^2+3}{4 \left(w^2+1\right)},\\
\\
x_9 &=-\frac{(v-1)^4 (8 v-5) (18 v-5)}{4 (3 v-2)^2},&c_9=- (-5 + 8 v) (-5 + 18 v),&    &v_9=\frac{5 \left(w^2+1\right)}{2 \left(4 w^2+9\right)}\\
\\
x_{10} &=\frac{1}{8} \left(-4 v^2+4 v+1\right),&c_{10}=-2(-1 - 4 v + 4 v^2),&   &v_{10}=\frac{w^2-4 w+2}{2 \left(w^2+2\right)}.\\
\end{align*}
}

\subsection{Second part of the list.}
\mbox{}
{\footnotesize
\begin{align*}
x_{11} &=-\frac{(v-1)^2 (2 v-5)^2 \left(36 v^2-70 v+25\right)}{(6 v-7)^2},&c_{11}=-25 + 70 v - 36 v^2, &   &v_{11}=\frac{w^2-6 w+34}{w^2+36},\\
\\
x_{12} &=-\frac{4 (v-1)^3 v (2 v+1)^2}{(2 v-3)^2},&c_{12}=-16 (-9 - 28 v + 28 v^2), &    &v_{12}=-\frac{2 (3 w-14)}{w^2+28},\\
\\
x_{13} &=-\frac{4 (v-1)^3 v (10 v-1)}{10 v-9},& c_{13}=(-9 + 10 v) (-1 + 10 v), &     &v_{13}=\frac{(3 w-1) (3 w+1)}{10 (w-1) (w+1)},\\
\\
x_{14} &=-\frac{27}{2} (v-1) v,&c_{14} =-6 (-3 + v) (2 + v), & &v_{14}=-\frac{2 (w-3) (w+3)}{w^2+6},\\
\\
x_{15} &=-\frac{1}{64} \left(16 v^2-16 v+1\right)^2,&c_{15}=(7 - 128 v + 128 v^2),     &  & v_{15}=-\frac{w^2+80 w+1152}{8 \left(w^2-128\right)},\\
\\
x_{16} &=\frac{(v-1)^2 (4 v-1)^2 (10 v-1)}{8 (3 v-1)},& c_{16}= 2  (-1 + 3 v) (-1 + 10 v) , &    &v_{16}=\frac{2 w^2-1}{2 \left(3 w^2-5\right)}.\\
\end{align*}
}

\subsection{ Third part of the list.}
\mbox{}

{
\footnotesize
\begin{align*}
x_{17} &=-\frac{(v-1)^2 (2 v+1)^2 (8 v+1)}{8 (v-3)},&c_{17}=  -2 (-3 + v) (1 + 8 v),&  &v_{17}=\frac{6 w^2-1}{2 \left(w^2+4\right)},\\
\\
x_{18} &=\frac{4 (4 v-3)^2 (10 v-9) \left(18 v^2-26 v+9\right)^2}{(6 v-5)^2 (18 v-13)},&c_{18}= (-9 + 10 v) (-13 + 18 v),&    &v_{18}=\frac{9 w^2-13}{2 \left(5 w^2-9\right)}.\\
\end{align*}
}

All of them have rank at least $1$,
and  all of them are non-isomorphic to each other.
That the rank of these three curves is at least $1$ over $\Q(w)$ can be proved using a specialization argument
since the specialization map is a homomorphism, see \cite[Theorem 11.4]{S}.
That all of them are non-isomorphic to each other can be proved by comparing their $j$-invariants.

\section{Rank $2$ families  for the torsion group  $\Z/8\Z$}\label{sec:rank two}

In this section and in the next one we present curves with torsion group  $\Z/8\Z$ and rank $2$ over  $Q(u)$.
The search for these curves is made again by looking for  solvable quadratic conditions  leading  to new points.
Most  of these conditions can be obtained also by comparing the $j$-invariants of the different curves with rank $1$, for example from $v_3(u)=v_{14}(s)$ we deduce  $w_3$.
But we have followed  mainly the first alternative in the presentation  because in this way we get directly the $x$-coordinates of the new point.

 \subsection{First Dujella and Peral case}
We present here some details for the family  in which we have found three subfamilies with torsion  $\Z/8\Z$ and
generic rank at least $2$. It corresponds to the third  entry in the table of rank $1$ families above. By inserting in the general family $y^2= x^3+  A_{8}(v) x^2+  B_{8}(v) x$ the value $v=v_3(w)$ and clearing denominators we get  the rank $1$ family given by $y^2= x^3+ a_{3}(w) x^2+  b_{3}(w) x$, where
\begin{align*}
a_{3}(w)&=-31 - 148 w^2 + 214 w^4 - 116 w^6 + 337 w^8,\\
b_{3}(w)&=256 (-1 + w)^4 (1 + w)^4 (1 + 3 w^2)^4.
\end{align*}
The $X$-coordinate of a  point of infinite order is
\[
X = -256 (-1 + w)^3 w^2 (1 + w)^3.
\]

By searching on several  homogeneous spaces of this curve, we have  found three  conditions that lead to new points.

We list  the values  $X_1$, $X_2$ and $X_3$ jointly with the specialization of the parameter that converts them into a point of the specialized  curve
\begin{align*}
X_1 &=\frac{(-1 + w)^2 (1 + w)^2 (5 + 7 w^2)^2 (11 + 25 w^2)}{16},&w_1=\frac{11-u^2}{10 u},\\
X_2 &=\frac{(-1 + w)^2 (1 + w)^2 (1 + 11 w^2)^2 (7 + 29 w^2)}{16 w^2},&w_2=\frac{29 - 12 u + u^2}{-29 + u^2},\\
X_3 &=-27 (w-1) (w+1) \left(w^2+3\right)^2 \left(3 w^2+1\right),&w_3=\frac{u^2-12 u+15}{u^2-15}.
\end{align*}

 In this way  we get  three curves of rank  $2$ over $\Q(u)$. The two  curves corresponding to $w_1$ and $w_2$ were reported in \cite{DP1} and \cite{DP2}. The curve corresponding to  $w_3$ is new.

 Once we insert $w_1$ into the coefficients $ a_{3}, b_{3}$, we get as new coefficients  $aa_{1}, bb_{1}$ given by
\begin{align*}
aa_{1} &=   337 u^{16} - 41256 u^{14}  + 4047356 u^{12}  -288332632 u^{10 } + 2363813190 u^8 \\
& - 34888248472 u^6 + 59257339196 u^4- 73087520616 u^2+72238942897, \\
bb_{1} &= 256\, (363 + 34 u^2 + 3 u^4)^4\,(11 + u)^4\,(-11 + u)^4\,(-1 + u)^4\, (1 + u)^4.
\end{align*}
The $X$-coordinates of  two  independent infinite  order points   are
\begin{align*}
 &-16 (u-11)^3 (u-1) (u+1)^3 (u+11) \left(3 u^4+34 u^2+363\right)^2, \\
&\frac{(-11 + u)^2 (-1 + u)^2 (1 + u)^2 (11 + u)^2 (11 + u^2)^2 (847 +
   346 u^2 + 7 u^4)^2}{64 u^2}.
\end{align*}

 \subsection{Second Dujella and Peral case.}
By inserting   $w_2$ into the coefficients $ a_{3}, b_{3}$, we get as new coefficients  $aa_{2}, bb_{2}$ given by
\begin{align*}
aa_{2} &=   500246412961 - 2069985157080 u + 3162080774436 u^2 -
 2895517882032 u^3 +\\& 1873181389706 u^4 - 906769167048 u^5 +
 333391978480 u^6 - 93284915496 u^7 +\\& 19860033555 u^8 -
 3216721224 u^9 + 396423280 u^{10} - 37179432 u^{11} +\\& 2648426 u^{12} -
 141168 u^{13} + 5316 u^{14} - 120 u^{15} + u^{16}, \\
bb_{2} &=256 (-6 + u)^4 u^4 (-29 + 6 u)^4 (841 - 522 u + 137 u^2 - 18 u^3 +
   u^4)^4.
\end{align*}

The $X$-coordinates of  two  independent infinite  order points   are
\begin{align*}
&16 (u-6) u (6 u-29)^3 \left(u^4-18 u^3+137 u^2-522 u+841\right)^2,  \\
&\frac{(-6 + u)^2 u^2 (-29 + 6 u)^2 (87 - 29 u + 3 u^2)^2 (2523 - 1914 u +
   541 u^2 - 66 u^3 + 3 u^4)^2}{4 (29 - 12 u + u^2)^2}.
\end{align*}

 \subsection{The first new curve.}
Now  we insert $w_3$ into the coefficients $ a_{3}, b_{3}$. We get as new coefficients  $aa_{3}, bb_{3}$ given by

 \begin{align*}
aa_{3} &= 2562890625 - 20503125000 u + 58638937500 u^2 - 98524350000 u^3 +\\&
 112751696250 u^4 - 92004903000 u^5 + 54062154000 u^6 -
 22880209320 u^7 +\\& 6966724707 u^8 - 1525347288 u^9 + 240276240 u^{10} -
 27260712 u^{11} +\\& 2227194 u^{12} - 129744 u^{13} + 5148 u^{14} -
 120 u^{15} + u^{16}, \\
bb_{3} &=20736 (-6 + u)^4 u^4 (-5 + 2 u)^4 (75 - 15 u + u^2)^4 (3 - 3 u +
   u^2)^4.
\end{align*}
The $X$-coordinates of  two  independent infinite  order points   are
\begin{align*}
 &1728 (-6 + u)^3 u^3 (-5 + 2 u)^3 (15 - 12 u + u^2)^2,  \\
&81 (-6 + u) u (-5 + 2 u) (75 - 15 u + u^2) (3 - 3 u + u^2) (225 -
   90 u + 21 u^2 - 6 u^3 + u^4)^2.
\end{align*}

 \subsection{The first MacLeod curve with rank $2$.}

 We give here some details for the curve  in which MacLeod  discovered  a
curve having rank  $2$. It corresponds to the entry  $v_{12}$ in the table for rank $1$ curves above. By inserting $v=v_{12}(w)$  in the general family $y^2= x^3+  A_{8}(v) x^2+  B_{8}(v) x$, we get  the rank $1$ curve given by $y^2= x^3+ a_{12}(w) x^2+  b_{12}(w) x$ where
\begin{align*}
a_{12}(w)&=614656 - 1053696 w + 363776 w^2 - 59136 w^3 - 7328 w^4 + 2112 w^5 \\
&\,\,\,\,\mbox{}+ 464 w^6 + 48 w^7 + w^8,\\
b_{12}(w)&=256 w^4 (6 + w)^4 (-14 + 3 w)^4.
\end{align*}
The $X$-coordinate of the  point of infinite order is
\[
X = -\frac{8 w^3 (w+6)^3 (3 w-14) \left(w^2-12 w+84\right)^2}{\left(3
   w^2+12 w+28\right)^2}.
\]

Now MacLeod imposed $-((256 w^3 (6+w)^2 (-14+3 w)^3)/(14+w)^2)$ as the $X$-coordinate of a new point. This is the same as  to spezialize  $w_4=-\frac{(u-28) (u+28)}{2 (u-63)}$.
Using this value into the coefficients $ a_{12}, b_{12}$, we get as new coefficients  $aa_{4}, bb_{4}$ given by
  \begin{align*}
aa_{4} &= 1058387660788345388204032 - 141209336315730168643584 u +\\&
 7118408590330053918720 u^2 + 46091099527055278080 u^3 - \\&
 20521534612217970688 u^4 + 473831305485189120 u^5 +
 19585996741025792 u^6 - \\&545185218600960 u^7 - 18026420955648 u^8 +
 234415749120 u^9 +22250170880 u^{10} - \\& 242597376 u^{11} -
 14269120 u^{12} + 276096 u^{13} + 1632 u^{14} - 96 u^{15} + u^{16}, \\
bb_{4} &=4096 (-63 + u)^4 (-28 + u)^4 (-14 + u)^4 (2 + u)^4 (28 + u)^4 (42 +
   u)^4 (-98 + 3 u)^4.
\end{align*}
The $X$-coordinates of  two  independent infinite  order points   are
\begin{align*}
 &16 (u-63) (u-28)^3 (u-14)^3 (u+2)^3 (u+28)^3 (u+42) (3 u-98)\times\\&\frac{
   \left(u^4+24 u^3-2744 u^2-61152 u+3133648\right)^2}{\left(3
   u^4-24 u^3-3080 u^2+4704 u+1103088\right)^2},  \\
&-\frac{1024 (u-63)^2 (u-28)^3 (u-14)^2 (u+2)^2 (u+28)^3 (u+42)^3 (3
   u-98)^3}{\left(u^2-28 u+980\right)^2}.
\end{align*}

 \subsection{The second MacLeod curve with rank $2$.}

 We present here some details for the curve  in which MacLeod \cite{ML}  discovered  a second
curve having rank  $2$. It corresponds to the entry  $v_{13}$ in the table for rank $1$ curves above. By inserting $v=v_{13}(w)$  in the general family $y^2= x^3+  A_{8}(v) x^2+  B_{8}(v) x$, we get  the rank $1$ curve given by $y^2= x^3+ a_{13}(w) x^2+  b_{13}(w) x$ where
\begin{align*}
a_{13}(w)=&2 (431 + 3524 w^2 - 3814 w^4 + 3524 w^6 + 431 w^8)\\
b_{13}(w)=&(-3 + w)^4 (3 + w)^4 (-1 + 3 w)^4 (1 + 3 w)^4.
\end{align*}
The $X$-coordinate of the  point of infinite order is
\[
X = (-3 + w)^3 w^2 (3 + w)^3 (-1 + 3 w) (1 + 3 w).
\]

Now MacLeod \cite{ML} imposed
\[
\frac{10 (w-3)^2 w^2 (w+3)^2 \left(13 w^2+3\right)^2}{17 w^2+7}
\]
as the $X$-coordinate of a new point. This is the same as  to spezialize  $w_5=-\frac{3 u^2-80 u+510}{u^2-170}$.
Using this value into the coefficients $ a_{13}, b_{13}$, we get as new coefficients  $aa_{5}, bb_{5}$ given by

  \begin{align*}
aa_{5} &= 4 u^{16}-768 u^{15}+68736 u^{14}-3816768 u^{13}+147831608
   u^{12}-\\& 4261407840 u^{11}+95281085176 u^{10}-1698380209632
   u^9+ 24531870965502 u^8-\\&288724635637440 u^7+2753623361586400
   u^6-20936296717920000 u^5+\\&123470437317680000 u^4-541926476217600000
   u^3+1659119942784000000 u^2-\\& 3151401008640000000
   u+2790302976400000000, \\
bb_{5} &=u^4 (3 u-40)^4 (4 u-51)^4 \left(u^2-24 u+136\right)^4 \left(2 u^2-60
   u+425\right)^4.
\end{align*}

The $X$-coordinates of  two  independent infinite  order points   are
\begin{align*}
 &-u (3 u-40) (4 u-51)^4 \left(2 u^2-60 u+425\right) \left(u^3-44 u^2+660 u-3400\right)^2,  \\
&-u (3 u-40) (4 u-51)^2 \left(2 u^2-60 u+425\right) \left(35 u^3-1236 u^2+14620
   u-57800\right)^2.
\end{align*}

 \subsection{The second new curve with rank $2$.}

 We give here some details for the curve  of rank $1$ in which we have found another  new
curve with rank  $2$.
 It corresponds to the entry  $v_{17}$ in the table for rank $1$ curves above. By inserting $v=v_{17}(w)$  in the general family $y^2= x^3+  A_{8}(v) x^2+  B_{8}(v) x$, we get  the rank $1$ curve given by $y^2= x^3+ a_{17}(w) x^2+  b_{17}(w) x$ where
\begin{align*}
a_{17}(w)=&2 (1169 - 3956 w^2 + 3704 w^4 - 2216 w^6 + 674 w^8)\\
b_{17}(w)=&(-3 + 2 w)^4 (3 + 2 w)^4 (-1 + 6 w^2)^4.
\end{align*}
The $X$-coordinate of the  point of infinite order is
\[
X =\frac{1}{4} w^2 (2 w-3)^2 (2 w+3)^2 \left(7 w^2+3\right)^2.
\]

 Observe that in order to have   \[X=-(27/2) (-3+2 w) (3+2 w) (4+w^2)^2 (-1+6 w^2)\] as a new point, we have to solve $30(3+2 w^2)$ equal to a rational square.
It is enough   to specialize to
\[
w_6=\frac{3 \left(u^2-20 u+60\right)}{2 \left(u^2-60\right)}.
\]
Now we insert $w_6$ in  the coefficients $ a_{17}, b_{17}$,  and we get as new coefficients  $aa_{6}, bb_{6}$ given by

  \begin{align*}
aa_{6} &= 625 u^{16}-180000 u^{15}+17872800 u^{14}-1010171520 u^{13}+\\&37753002432
   u^{12}- 973296787968 u^{11}+17592030254592 u^{10}-\\&225415897049088
   u^9+ 2063161668920832 u^8-13524953822945280 u^7+\\& 63331308916531200
   u^6- 210232106201088000 u^5+489278911518720000
   u^4-\\&785509373952000000 u^3+ 833873356800000000
   u^2\\& -503884800000000000 u+104976000000000000, \\
bb_{6} &=6879707136 (u-10)^4 (u-6)^4 u^4 \left(u^2-36 u+300\right)^4 \left(5
   u^2-36 u+60\right)^4.
\end{align*}

The $X$-coordinates of  two  independent infinite  order points   are
\begin{align*}
&104976 (u-10)^2 (u-6)^2 u^2 \left(u^2-20 u+60\right)^2\times \\&
 \frac{ \left(5
   u^4-168 u^3+2088 u^2-10080 u+18000\right)^2}{\left(u^2-60\right)^2},  \\
&1944 (u-10) (u-6) u \left(u^2-36 u+300\right) \left(5 u^2-36
   u+60\right) \times\\& \left(5 u^4-72 u^3+552 u^2-4320 u+18000\right)^2.
\end{align*}

 \subsection{The third new curve with rank $2$.}

 We describe  here some details for the curve  of rank $1$ in which we have found another  new
curve with rank  $2$.
 It corresponds to the entry  $v_{18}$ in the table for rank $1$ curves above. By inserting $v=v_{18}(w)$  in the general family $y^2= x^3+  A_{8}(v) x^2+  B_{8}(v) x$, we get  the rank $1$ curve given by $y^2= x^3+ a_{18}(w) x^2+  b_{18}(w) x$ where
\begin{align*}
a_{18}(w)=&2 (-3713 + 5492 w^2 - 1462 w^4 - 1004 w^6 + 431 w^8)\\
b_{18}(w)=&((-5 + w^2)^4 (-13 + 9 w^2)^4.
\end{align*}
The $X$-coordinate of the  point of infinite order is
\[
X =\frac{\left(3 w^2+1\right)^2 \left(9 w^4+70 w^2-63\right)^2}{w^2 \left(w^2+3\right)^2}.
\]
If we impose
\[
X=(4 (1+3 w^2)^2 (-13+9 w^2)^2)/(-3+7 w^2)
\]
as a new point, we have to solve $(-3+7 w^2)$ equal to a rational square.
This is done with  \[
w_7=\frac{u^2-6 u+21}{u^2-14 u+21}
\]
Now we use $w_7$ in  the coefficients $ a_{18}, b_{18}$  and we get as new coefficients  $aa_{7}, bb_{7}$ given by

{\tiny
  \begin{align*}
aa_{7} &=-2 u^{16}+384 u^{15}-30128 u^{14}+1278592 u^{13}-32804472 u^{12}+545481088 u^{11}-\\& 6133914960 u^{10}+47788256896
   u^9-261061974220 u^8+1003553394816 u^7-\\& 2705056497360 u^6+5051700355968 u^5-6379846519032 u^4+\\& 5221898865792
   u^3-2583961693488 u^2+691617999744 u-75645718722, \\
bb_{7} &=\left(u^2-56 u+147\right)^4 \left(u^2-8 u+3\right)^4 \left(u^4-32 u^3+278 u^2-672 u+441\right)^4.
\end{align*}
}

The $X$-coordinates of  two  independent infinite  order points   are
{\tiny
\begin{align*}
&\left(u^2-8 u+3\right)^2 \left(u^4-32 u^3+278 u^2-672 u+441\right) \left(2 u^5-55 u^4+508
   u^3-1834 u^2+3234 u-3087\right)^2,  \\
   &\frac{u^2 \left(u^2-56 u+147\right)^2 \left(u^2-8 u+3\right)^4 \left(u^4-32 u^3+278
   u^2-672 u+441\right)^3}{\left(u^5-22 u^4+262 u^3-1524 u^2+3465 u-2646\right)^2}.
\end{align*}
}

\subsection{Rank $2$ results}

\begin{theo}
The seven curves corresponding to the specializations $w_i$, $i=1,...,7$
have rank $2$ over $\Q(u)$, and the points listed in each case, jointly with the torsion points,  are generators for the full Mordell-Weil group.
\end{theo}

First we have proved that the curves have  rank at least $2$ using the specialization
theorem \cite[Theorem 11.4]{S}.
Then we have used the Gusi\'c-Tadi\'c algorithm \cite[Theorem 1.3]{GT}
to find an injective specialization,  and mwrank and magma to compute the rank and the Mordell-Weil group
for the specialized curves. The specializations listed below prove that the rank is exactly $2$
and that the points listed for each curve generate, jointly with the torsion points, the full Mordel-Weil group.
In some cases, the Gusi\'c-Tadi\'c algorithm indicated that the originally found points generated
a subgroup of order $2$, so by "halving" corresponding points, we found the generators of the full group.

\begin{align*}
 \hbox{$w_i$}&  &\hbox{$u$ value}\\
1& &22\\
2 & & 19\\
3& &11\\
4& &17\\
5 & &3\\
6& &-48\\
7& &10\\
\end{align*}

\section{Infinite families of rank $3$ for the  torsion group $\Z/8\Z$}

\subsection{First example.}\label{section:first example}
It can be proved that there exist infinitely many elliptic curves with torsion group $\Z/ 8 \Z $ parametrized by the points on an elliptic curve with positive rank.
For example,  it  is enough to see that the equation $w_1(r)=w_2(s)$, i.e.,
 \[
\frac{11- r^2}{10 r}=\frac{29 - 12 s + s^2}{-29 + s^2}
\]
  has infinitely many solutions. This is the same as to solve
  \[
r^2 s^2-29 r^2+10 r s^2-120 r s+290 r-11 s^2+319=0
  \]
   in rationals, so   the discriminant  $\Delta =3509 + 62 r^2 + 29 r^4$ has to be a square.

   Observe that this is the same as imposing
   \begin{align*}
   \frac{(u-11)^2 (u-1)^2 (u+1)^2 (u+11)^2 \left(11 u^4-142 u^2+1331\right)^2 \left(29 u^4+62
   u^2+3509\right)}{16 \left(u^2-11\right)^2}
    \end{align*}
   as a new point on the curve of rank $2$ corresponding to $w_1$.

   But  $t^2=3509 + 62 r^2 + 29 r^4$ has a solution, for example $(r,t)=(1,60)$, hence it is equivalent to the cubic $Y^2=X^3- 463 X^2 +45936 X$ whose rank  is $2$ as proved with mwrank \cite{Cr}.
   This, jointly with the independence of the corresponding points, implies the  existence of infinitely many solutions parametrized by the points of the elliptic curve,
   see \cite{Le1} or \cite{Ra} for this kind of argument.

   \subsection{Other examples.}\label{sec:other examples}
For the rank $2$ curve corresponding to $w_1$, imposing
 \[
X=-27 (u-11) (u-1) (u+1) (u+11) \left(u^4+278 u^2+121\right)^2 \left(3 u^4+34 u^2+363\right)\]
as a new point it is equivalent to solve
 \[
 15 (121 + 118 u^2 + u^4) =t^2
 \]
This equation  has for  example  the solution $(u,t)(=1,60)$.
So it is equivalent to the elliptic curve
$Y^2 = X^3 + 1770X^2 -108900X-192753000$ whose rank is $1$.


For the rank $2$ curve corresponding to $w_2$, imposing

{\tiny \[
X=27 (u-6) u (6 u-29) \left(u^4-18 u^3+137 u^2-522 u+841\right) \left(u^4-6 u^3+7 u^2-174 u+841\right)^2
\]
}
as a new point is equivalent to solve
 \[
 3(2523-870 u+151 u^2-30 u^3+3 u^4) =t^2,
 \]
 having the solution $(u,t)=(0,87)$, so the
  quartic can be shown equivalent to the elliptic curve
  $Y^2 = X^3 +453X^2 -37584X-817452$ whose rank is $2$.

  \section{Examples of curves with high rank}

         The highest known rank of an elliptic curve over $\Q$ with  torsion group $\Z/8\Z$ is  $6$. The first was discovered by Elkies  in 2006,  the second was found by Dujella MacLeod and Peral in 2013, and the third  and fourth by Voznyy in 2021.

         The second curve with rank $6$ corresponds to $w=-\frac{261}{70}$ in the curve number $12$.
         The third curve with rank  $6$ corresponds to $w=\frac{1327}{989}$ in the curve number $13$.
         See \cite{D} for the details of these curves.

The following list includes examples of rank $5$ curves  found in the rank $1$ curves.
First column indicates the number of the  curve, and the second the value(s) of the parameter producing a rank $5$ curve.
The details on the curves, including the authors and years of discoveries,
can be found on the web page \cite{D}.

\begin{align*}
 \hbox{Curve number}&  &\hbox{$w$ values}\\
1 \hspace{1cm}&& \frac{72}{19},\, \frac{101}{145},\\
2 \hspace{1cm}&& \frac{317}{10},\, \frac{235}{46}, \, \frac{309}{130},\\
3 \hspace{1cm}&& \frac{73}{83}, \, \frac{37}{157},\, \frac{131}{419},\, \frac{699}{1226},\, \frac{166}{121},\\
4 \hspace{1cm}&& \frac{245}{12},\,\frac{95}{396},\, \frac{-87}{28},\\
6 \hspace{1cm}&& \frac{100}{29},\, \frac{-28}{79},\, \frac{304}{55},\\
7 \hspace{1cm}&& \frac{79}{431},\\
9 \hspace{1cm}&& \frac{287}{109},\, \frac{65}{71},\\
10 \hspace{1cm}&& \frac{21}{95},\, \frac{195}{154},\\
11 \hspace{1cm}&& \frac{94}{31},\, \frac{103}{136}, \, \frac{508}{201},\\
12 \hspace{1cm}&& \frac{3}{22},\, \frac{117}{40}, \,  \frac{48}{209}, \, \frac{24}{43},  \, \frac{-153}{4},  \, \frac{193}{2}, \, \frac{533}{126},  \, \frac{1440}{319},  \, \frac{-7446}{2773},\, \frac{2365}{426},  \\
13 \hspace{1cm}&&\frac{501}{2407},\,  \frac{77}{188},\,  \frac{427}{1341},\\
14 \hspace{1cm}&& \frac{838}{331},\,   \frac{484}{283},\,   \frac{1538}{631},\,   \frac{382}{121},\,   \frac{367}{94},\,   \frac{1226}{477},\\
16 \hspace{1cm}&&\frac{657}{262}, \,\frac{283}{82},\,\frac{4969}{2796},\,\frac{2742}{839},\, \frac{906}{719},\, \frac{762}{521}, \\
17 \hspace{1cm}&& \frac{1557}{1538},\,  \frac{27}{382}, \,\frac{489}{670},\,\frac{811}{1351},\,\frac{198}{97} \\
18 \hspace{1cm}&& \frac{89}{51},\, \frac{1099}{1371}.\\
\end{align*}


\section{Rank $1$   curves for  torsion group $\Z/2\Z \times  \Z/6\Z $.}

We have found over forty curves having rank $1$ for this torsion group,
but we list only those in which we have found a further quadratic sections producing
the five non-isomorphic rank  $2$ curves.

We present in the next section five non-isomorphic curves of rank $2$, the four already known
and another one that is new.
The first four curves  of rank $2$ for this torsion group  were discovered using Diophantine triples or a model due to Hadano \cite{DP1}, \cite{DP2}. The fifth one has been found by using two consecutive quadratic sections.    As said before, the seven cases presented by MacLeod \cite{ML}, are all isomorphic to one of the first three curves mentioned.

In this presentation, we  get all these curves by imposing pair of consecutive quadratic  sections. The data for the rank $1$ curves in which we found rank $ 2$ curves are in the next list. In the first rank 1 curve,
we have found two rank $2$ curves.

{\footnotesize
\begin{align*}
x_{1} &=16 (-2 + v) (1 + v)^2, &c_1= 3 (-2 + v) v,& &v_{1}=-\frac{6}{w^2-3},\\
x_{2} &=\frac{64 (v-1)^2 (v+1)^3}{(v+5)^2},&c_2=-6 (-7 + v) (3 + v),   &  &v_{2}=\frac{3 \left(14 w^2-1\right)}{6 w^2+1},\\
x_{3} &=(1 + v)^2 (7 - 4 v + v^2)^2,&c_3=6 - 2 v + v^2,& &v_{3}=\frac{2 w^2-4 w-1}{2 w-3},\\
x_{4} &=\frac{64}{3} (v-1) (v+1) (3 v-1)^2,&c_4= 6 (-1 + 2 v)  (-1 + 7 v) ,&  &v_{4}=\frac{6 w^2-1}{12 w^2-7}.\\
\\
\end{align*}
}


\section{Rank $2$   curves for  torsion group $\Z/2\Z \times  \Z/6\Z $.}

\subsection{First Dujella-Peral curve with rank $2$.}\label{sec:8.1}
As said  before, this curve was first discovered using Diophantine triples \cite{DP1}, \cite{DP2}. Here we  give the details  using  a pair of quadratic sections. The curve of rank $1$  corresponding to $v_1$ in the previous table is the following:

\begin{align*}
a_{1} &= -4779 - 4644 w^2 + 1134 w^4 + 60 w^6 + 37 w^8,\\
b_{1} &= 32 (-3 + w)^3 (3 + w)^3 (-3 + w^2) (3 + w^2)^3 (3 + 5 w^2).
\end{align*}
A point of infinite order has the following $X$-coordinate
\begin{align*}
X &=-32 (-3 + w)^2 w^2 (3 + w)^2 (-3 + w^2).\\
\end{align*}
Now we impose $-(-3+w)^2 (3+w)^2 (-3+w^2) (9+7 w^2)$ as a new point, this is the same as specializing to $w_1=\frac{3 \left(u^2-8 u+14\right)}{u^2-14}$. From here, we get the rank $2$ curve whose main data are:

 {\tiny
  \begin{align*}
aa_{1} &= 1475789056 - 6324810240 u + 12303261824 u^2 - 14934296832 u^3 +
 12836014912 u^4 - 8279778528 u^5 + \\&  4113507272 u^6 - 1590783936 u^7 +
 480725533 u^8 - 113627424 u^9 + 20987282 u^{10} - 3017412 u^{11} +
 334132 u^{12} -  \\& 27768 u^{13} + 1634 u^{14} - 60 u^{15} + u^{16}, \\
   bb_{1} &=-27 (-4 + u)^3 u^3 (-7 + 2 u)^3 (196 - 336 u + 152 u^2 - 24 u^3 +
   u^4) (196 - 168 u + 62 u^2 - 12 u^3 + u^4)^3\\& (392 - 420 u +
   169 u^2 - 30 u^3 + 2 u^4).
\end{align*}
}

The $X$-coordinates of  two  independent infinite  order points   are
{\tiny
\begin{align*}
&-27 (-4 + u)^2 u^2 (-7 + 2 u)^2 (14 - 8 u + u^2)^2 (196 - 336 u +
   152 u^2 - 24 u^3 + u^4),  \\
   &-\frac{27}{4} (u-4)^2 u^2 (2 u-7)^2 \left(u^2-7 u+14\right)^2 \left(u^4-24 u^3+152 u^2-336 u+196\right).
\end{align*}
}


\subsection{Second  Dujella-Peral curve with rank $2$.}

Now we impose \[(-3+w) (3+w) (-3+w^2) (9+7 w^2)^2\] in the rank $1$ curve for $v_1$.
This is the same as specializing to  $w_2= \frac{u^2-8 u+6}{u^2-6}$. We get the rank $2$ curve whose coefficients  are:

 {\tiny
  \begin{align*}
aa_{2} &=-3359232 + 2239488 u + 6905088 u^2 - 11695104 u^3 + 6925824 u^4 -
 2494368 u^5 + 3007512 u^6 - \\& 3509088 u^7 + 2015437 u^8 - 584848 u^9 +
 83542 u^{10} - 11548 u^{11} + 5344 u^{12} - 1504 u^{13} + 148 u^{14} +
 8 u^{15} - 2 u^{16}, \\
   bb_{2} &=(-3 + u)^3 (-2 + u)^3 (1 + u)^3 (6 + u)^3 (36 - 60 u + 43 u^2 -
   10 u^3 + u^4) (36 - 24 u + 10 u^2 - 4 u^3 + u^4)^3 \\& 36 + 48 u -
   56 u^2 + 8 u^3 + u^4).
\end{align*}
}

The $X$-coordinates of  two  independent infinite  order points   are
{\tiny
\begin{align*}
&(-3 + u)^2 (-2 + u)^2 (1 + u)^2 (6 + u)^2 (6 - 8 u + u^2)^2 (36 +
   48 u - 56 u^2 + 8 u^3 + u^4),  \\
   &\frac{1}{4} (u-3) (u-2) (u+1) (u+6) \left(u^4+8 u^3-56 u^2+48 u+36\right) \left(2 u^4-14 u^3+53 u^2-84 u+72\right)^2.
   \end{align*}
}


\subsection{Third Dujella-Peral curve with rank $2$.}\label{sec:8.3}

We use $v_2$ and we  get the rank $1$ curve whose coefficients are:
\begin{align*}
a_{2} &= 121 - 2136 w^2 - 5184 w^4 + 273024 w^6 - 1223424 w^8,\\
b_{2} &= 128 (-1 + 3 w)^3 (1 + 3 w)^3 (1 + 6 w^2) (-1 + 24 w^2)^3 (-7 + 48 w^2).
\end{align*}
A point of infinite order has the following $X$-coordinate
\begin{align*}
X &=\frac{128 (3 w-1)^2 (3 w+1)^2 \left(6 w^2+1\right) \left(24 w^2-1\right)^3}{\left(36 w^2+1\right)^2}.\\
\end{align*}
Now we impose $4 (-1 + 3 w)^2 (1 + 3 w)^2 (1 + 36 w^2) (-7 + 48 w^2)$ as $X$-coordinate for a new point. This is the same as using $w_3=\frac{u^2-30 u+180}{3 \left(u^2-180\right)}$ in the preceding rank $1$ curve. We get the rank $2$ curve whose coefficients are:

 {\tiny
  \begin{align*}
aa_{3} &=1101996057600000000 - 587731230720000000 u - 901187887104000000 u^2 +
 2262765238272000000 u^3 -\\& 1225425058007040000 u^4 +
 335908714991616000 u^5 - 57791877967872000 u^6 +
 6886398457405440 u^7 - \\& 597067735693824 u^8 + 38257769207808 u^9 -
 17836999372{13 }- 26496 u^{14}- 96 u^{15}+ u^{16},  \\
   bb_{3} &=5971968 (-15 + u)^3 (-12 + u)^3 u^3 (32400 - 17280 u + 2232 u^2 -
   96 u^3 + u^4)^3\\& (32400 - 4320 u + 288 u^2 - 24 u^3 + u^4) (32400 +
   34560 u - 5544 u^2 + 192 u^3 + u^4).
\end{align*}
}

The $X$-coordinates of  two  independent infinite  order points   are
{\tiny
\begin{align*}
&\frac{18432 (u-15)^2 (u-12)^2 u^2 \left(u^4-96 u^3+2232 u^2-17280 u+32400\right)^3 \left(u^4-24 u^3+288 u^2-4320
   u+32400\right)}{\left(u^2-24 u+180\right)^4}, \\
   &-15552 (u-15)^2 (u-12)^2 u^2 \left(u^2-24 u+180\right)^2 \left(u^4+192 u^3-5544 u^2+34560 u+32400\right).
   \end{align*}
}


\subsection{Fourth Dujella-Peral curve with rank $2$.}

We use $v_3$ and we  get the rank $1$ curve whose coefficients are:
\begin{align*}
a_{3} &= 96 - 480 w + 1584 w^2 - 3084 w^3 + 3001 w^4 - 1440 w^5 + 306 w^6 -
 12 w^7 - 3 w^8,\\
b_{3} &=16 (-3 + w) (-2 + w)^3 (1 + w)^3 (-3 + 2 w) (-2 + 3 w) (1 - 3 w +
   w^2)^3.
\end{align*}
A point $P$ of infinite order has the following $X$-coordinate
\begin{align*}
X &=\frac{4 (w-2)^2 (w+1)^2 \left(w^4-8 w^3+24 w^2-29 w+13\right)^2}{(2 w-3)^2}.\\
\end{align*}
Now we impose $(-2+w)^3 (1+w) (-3+2 w) (-2+3 w) (1-3 w+w^2) $ as $X$-coordinate for a new point. This is the same
as using $w_4=-\frac{4 u+9}{u^2-3}$ in the preceding rank $1$ curve. We get the rank $2$ curve whose coefficients are:

 {\tiny
  \begin{align*}
aa_{4} &=-314928 - 7978176 u - 47134224 u^2 - 141974208 u^3 - 263196864 u^4 -
 321113808 u^5 - 259493652 u^6 - \\& 128609568 u^7 - 23353995 u^8 +
 16908960 u^9 + 16006092 u^{10} + 6735888 u^{11} + 1706128 u^{12} +
 271104 u^{13} +\\&  27360 u^{14} + 1920 u^{15} + 96 u^{16},\\
   bb_{4} &=16 (u-6)^3 u (u+2)^3 (3 u+4) \left(u^2-3\right) \left(u^2+3 u+1\right)^3 \left(u^2+9 u+9\right)^3 \left(2 u^2+4
   u+3\right)^3\\&  \left(2 u^2+12 u+21\right) \left(3 u^2+8 u+9\right).
\end{align*}
}

The $X$-coordinates of  two  independent infinite  order points   are
{\tiny
\begin{align*}
&-4 (u-6) u (u+2) (3 u+4) \left(u^2+3 u+1\right) \left(u^2+9 u+9\right) \left(2 u^4+8 u^3+22
   u^2+48 u+45\right)^2 ,
\\
   &(6 - u) (2 + u) (1 + 3 u + u^2) (9 + 9 u + u^2) (3 + 4 u +
   2 u^2)^3 (21 + 12 u + 2 u^2) (9 + 8 u + 3 u^2).
   \end{align*}
}%
Note that the point $P$ is in $2E(\mathbb{Q})$, so we replaced it by a point $Q$ such that
$2Q = P \pmod{E(\mathbb{Q})_{\rm{tors}}}$ in order to obtain generators of the Mordell-Weil group.


\subsection{New curve with rank $2$.}

We use $v_4$ and we  get the rank $1$ curve whose coefficients are:

\begin{align*}
a_{4} &= 4048 - 22512 w^2 + 49248 w^4 - 50652 w^6 + 20493 w^8,\\
b_{4} &=432 (-1 + w)^3 (1 + w)^3 (-2 + 3 w)^3 (2 + 3 w)^3 (-7 +
   12 w^2) (-16 + 21 w^2).
\end{align*}
A point of infinite order has the following $X$-coordinate
\begin{align*}
X &=-64 (-1 + w) (1 + w) (-2 + 3 w) (2 + 3 w) (2 + 3 w^2)^2.\\
\end{align*}
Now we impose $-\frac{243}{4} (w-1)^2 (w+1)^2 \left(12 w^2-7\right) \left(21 w^2-16\right)$ as $X$-coordinate for a new point. This is the same as to use $w_5=\frac{4 \left(u^2+1\right)}{5 \left(u^2-1\right)}$ in the preceding rank $1$ curve. We get the rank $2$ curve whose coefficients are:

 {\tiny
  \begin{align*}
aa_{5} &=-675347 - 8801576 u^2 + 443877484 u^4 - 944081432 u^6 +
 22507829710 u^8 - 944081432 u^{10} + 443877484 u^{12}-\\&   8801576 u^{14} -
 675347 u^{16},
 \\
   bb_{5} &=6912 (-3 + u)^3 (3 + u)^3 (-1 + 3 u)^3 (1 + 3 u)^3 (11 + u^2)^3 (1 -
   5 u + u^2) (1 + 5 u + u^2)\times\\&  (1 + 11 u^2)^3 (17 + 734 u^2 + 17 u^4).
\end{align*}
}

The $X$-coordinates of  two  independent infinite  order points   are
{\tiny
\begin{align*}
&192 (u-3)^2 (u+3) (3 u+1) \left(u^2-5 u+1\right) \left(11 u^2+1\right)^2 \left(u^3+28 u^2+11
   u+8\right)^2,
   \\&
   243 (-3 + u)^2 (3 + u)^2 (-1 + 3 u)^2 (1 + 3 u)^2 (1 - 5 u +
   u^2) (1 + 5 u + u^2) (17 + 734 u^2 + 17 u^4).
   \end{align*}
}

\subsection{Rank $2$ results}

\begin{theo}
The five curves corresponding to the specializations $w_i$, $i=1,...,5$ have rank $2$ over $\Q(u)$,
and the points listed in each case, jointly with the torsion points, are generators for the full Mordell-Weil group.
\end{theo}

As for the torsion group $\Z/8\Z$,
we use the Gusi\'c-Tadi\'c algorithm again to find injective specializations and mwrank and magma
to compute the rank and the generators for the specialized curves. However, since for the curves with
the torsion group $\Z/ 2 \Z \times \Z/ 6 \Z$, the cubic polynomial factorizes into linear factors,
here we use \cite[Theorem 1.1]{GT}.
The specializations listed below prove that the rank is exactly $2$
and that the points listed for each curve generate, jointly with the torsion points, the full Mordell-Weil group.

\begin{align*}
 \hbox{$w_i$}&  &\hbox{$u$ value}\\
1& &15\\
2 & & 17\\
3& &22\\
4& &19\\
5 & &20\\
\end{align*}




\section{Infinite families of rank $3$ for the  torsion group $\Z/2\Z\times \Z/6\Z$} \label{sec:z26}

Imposing

\[
-\frac{27}{4} (u-4) u (2 u-7) \left(u^2-7 u+14\right)^4 \left(u^4-24 u^3+152 u^2-336
   u+196\right)
   \]
   as the $X$-coordinate of a a new point in the first  rank $2 $ curve  is equivalent to solve
   \[
   784 - 756 u + 293 u^2 - 54 u^3 + 4 u^4= t^2.
   \]
    This quartic has a rational point $ (u, t)=(0,28)$.
    The condition is equivalent to an elliptic curve with rank $1$.

In the same curve,  we have that imposing

{\tiny
\[(-4 + u) u (-7 + 2 u) (196 - 336 u + 152 u^2 - 24 u^3 + u^4) (196 -
   168 u + 62 u^2 - 12 u^3 + u^4) (392 - 420 u + 169 u^2 - 30 u^3 +
   2 u^4) \]
   }

   \noindent as the $X$-coordinate of a new point is the same as solving $196 - 420 u + 197 u^2 - 30 u^3 + u^4=t^2$,  which has a rational solution $(u,t)=(0,14)$ and can be seen to be equivalent to an elliptic curve with rank $1$.

   In the same curve, imposing
   \[
   \frac{27}{4} (u-4) u (2 u-7) \left(u^2-7 u+14\right)^2 \left(u^4-12 u^3+62 u^2-168 u+196\right)^2
   \]
    as the $X$-coordinate of a new point is the same that solving
     \[
    784 - 924 u + 383 u^2 - 66 u^3 + 4 u^4= t^2,
    \]
    which has a rational solution $(u, t)=(0, 28)$, and can be seen to be equivalent to an elliptic curve with rank $2$.

In the third curve, imposing

   \[
  -62208 (u-15)^2 (u-12)^2 u^2 \left(u^4-96 u^3+2232 u^2-17280 u+32400\right)^2
    \]
     as the $X$-coordinate of a new point is the same as solving
     \[
   32400 + 60480 u - 9432 u^2 + 336 u^3 + u^4=t^2
     \]
      which has a rational solution $(u, t)=(0, 180)$ and can be seen to be equivalent to an elliptic curve with rank $2$.

 \section{Examples of curves with high rank}

         The highest known rank of an elliptic curve over $\Q$ with  torsion group $\Z/2\Z\times \Z/ 6\Z$ is  $6$. The first was discovered by Elkies  in 2006,  the second was found by Dujella, Peral  and Tadi\'c in 2015 and the third  and fourth by Dujella and Peral  in 2020.
         The Elkies curve  with rank $6$ corresponds to $u=-\frac{5}{6}$ in the rank $2$ curve number $1$ and to $u=\frac{3}{4}$ in the rank $2$  curve number $3$.

                The following list includes examples of rank $5$ curves  found in the rank $2$ curves. First column indicates the number of the  curve, and the second the value(s) of the parameter producing a rank $5$ curve.
 See \cite{D} for the details of these curves.

\begin{align*}
 \hbox{Curve number}&  &\hbox{$v$ values}\\
1 \hspace{1cm}&&-\frac{5}{2},\\
2 \hspace{1cm}&&7, \, - 66,\, \frac{21}{17}, \,\frac{14}{9 },\, \frac{65}{27},  \\
4 \hspace{1cm}&& \frac{2}{5}, \, \frac{35}{4},\, -\frac{1}{10},\, -\frac{9}{62},\, \\
5 \hspace{1cm}&& \frac{13}{7},\, \frac{77}{6}.\\
\end{align*}





\section{On the rank: conjectures and heuristics.}\label{sec:heur}

There is an old conjecture predicting the existence of elliptic curves of arbitrarily high rank over $\Q$.
This has also been conjectured for each torsion group in the Mazur theorem.

But recently, some heuristic  predicts the existence of a universal bound for the rank of   the elliptic curves over $\Q$. See \cite{PPVW} for details. In that paper, the authors  state heuristic bounds  for each torsion group. In the case of torsion groups $\Z/ 8 \Z $ and $\Z/ 2 \Z \times \Z/ 6 \Z $,
their heuristic  claims  that $3$ is this bound, meaning  that only a finite number of elliptic curves with this torsion group would have rank over $\Q$ greater or equal to $4$.

As we said, we have several families of elliptic curves with torsion groups  $\Z/ 8 \Z$
and $\Z/ 2 \Z \times \Z/ 6 \Z $ and rank at least $3$ over $\Q$ parametrized by elliptic curves of positive rank.

Now observe that if we assume both the parity conjecture and  the heuristic bound, then as a consequence,  only a finite number of these curves can have root number equal to $1$, and so the root number would not be evenly  distributed over  these families. The numerical data that we have collected suggests that this might not be the case.

\begin{figure}[htp]
	\centering
	\includegraphics[scale = 0.58]{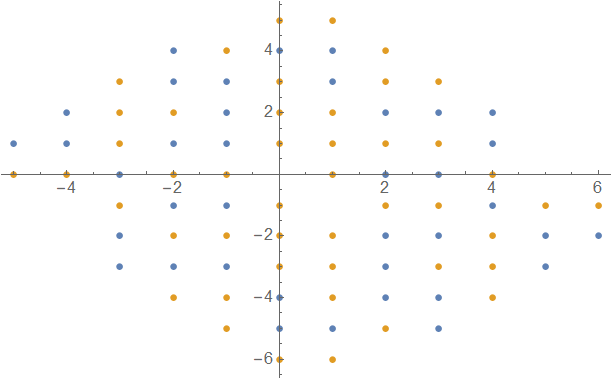}
	\caption{The blue (orange) point with coordinates $(n,m)$ represents the elliptic curve with root number one (minus one) that corresponds to the point $n P_1 + m P_2$ from the first example.}
	\label{fig:a}
\end{figure}

In our first example, see Section \ref{section:first example}, the elliptic curves with torsion group  $\Z/ 8 \Z$ and rank at least $3$ over $\Q$  are parametrized by rational points on the curve $$C:r^2 s^2-29 r^2+10 r s^2-120 r s+290 r-11 s^2+319=0.$$ This is a genus one curve (since it has two singular points with multiplicity $2$, $[1:0:0]$ and $[0:1:0]$, in its projective closure) with a rational point $[-1,0]$, so it is birationally equivalent to the elliptic curve in Weierstrass form $$E:Y^2 + X Y + Y = X^3 + X^2 - 1595 X  - 4768$$ by birational map $f$, which maps $[-1,0]$ to the point at infinity $\mathcal{O}\in E(\overline{\Q})$. This map is unique up to the composition with $[-1]$. We choose one such $f$. Elliptic curve $E$ has rank $2$ and a rational $2$-torsion. Denote the generators of the free part of the Mordell-Weil group by $P_1=(-57/4, 1043/8)$ and  $P_2=(42, -89)$ while the generators of the torsion are $T_1=(-3,1)$ and $T_2=(-39,19)$. The root numbers of elliptic curves corresponding (via $f$) to the points $n P_1+m P_2$ for small $n,m \in \Z$ are presented in Figure \ref{fig:a} (for the other choice of $f$, the figure would be centrally symmetric with respect to the origin). As we will show below, the points which differ by the point of order two correspond to the isomorphic elliptic curves, so they are not included in the figure.
There are $38$ curves with the root number $1$, and $46$ curves with the root number $-1$, which suggests that the root numbers are evenly distributed in the family.

The bottleneck of the root number computation is the factorization of the discriminant, which we need to determine the primes of bad reduction. One can speed up the computation by factoring the discriminant of the rank two elliptic curve over $\Q(w)$ from Section \ref{sec:rank two} (before specializing $w$) and using special number field sieve algorithm for factorization, but already for small values of $(n,m)$ one ends up factoring numbers with hundreds of digits. We are grateful to the members of Mersenne Forum (https://mersenneforum.org/) for their help with factorization that permitted us to compute the parity for the much larger number of curves than we initially thought it is possible (especially in our third example).

Moreover, note that the figure has a symmetry $(n,m) \leftrightarrow (-n+1,-m-1)$. This symmetry comes from the properties of two involutions $\tau_1$ and $\tau_2$ on $C$, which are birational automorphisms on $C$ that act on generic point $(r,s)$ as $\tau_1(r,s)=(r,\tilde{s})$ and $\tau_2(r,s)=(\tilde{r},s)$ where $\tilde{s}$ (respectively $\tilde{r}$) is the second solution of the quadratic equation in $s$ (respectively $r$) defined by $C$, i.e. $\tilde{s}=\frac{120 r}{r^2+10r-11}-s$. Since these functions extend to the regular involutions on the blowup of $C$ (which we identify with $E$ - the points in the blowup above the singular point $[0:1:0]$ are both rational and correspond to the points $P_1-P_2$ and $T_2$ in $E(\Q)$), we can identify them with involutions on $E$, which in general are of the form $T\mapsto -T +P$ for any point $P \in E(\overline{\Q})$ or of the form $T \mapsto T+P$ for $P$ of order $2$. Since both $\tau_1$ and $\tau_2$ have fixed points, they are of the form $T \mapsto -T+S_1$ and $T\mapsto -T+S_2$ for some points $S_1$ and $S_2$, respectively. One can check that, with our choice of generators, $S_1=P_1-P_2+T_1+T_2$ and $S_2=P_1-P_2$, thus the composition $\tau_1 \circ \tau_2$ is given by $T \mapsto T_1+T_2$. Since $C$ is equivalent to the condition $w_1(r)=w_2(s)$ (see Section \ref{sec:rank two}), it follows that the points in the same orbit under the action of group generated by $\tau_1$ and $\tau_2$ correspond to the isomorphic elliptic curve (with $\Z/8\Z$ torsion and rank at least $3$). In particular, it means that the points $T$ and $T+T_1+T_2$ correspond to the same elliptic curve, and likewise,  points $T$ and $-T+P_1-P_2$, which explains observed symmetry.

It remains to prove that points $T$, $T+T_1$ and $T+T_2$ correspond to the same curve. We first note that elliptic curves from Section \ref{sec:rank two} corresponding to parameters $w_1(r)$ and $w_1(-r)=-w_1(r)$ are the same. Moreover, we have another pair of involutions on C mapping $(r,s)\mapsto (-r,s')$ (since discriminant of the defining equation of $C$ with respect to $s$ is even function in $r$). More precisely, define $\psi_1(r,s)=(-r,\frac{s(r^2+10r-11)-120r}{r^2-10r-11})$ and $\psi_2(r,s)=\tau_1(\psi_1(r,s))$. Same as before, we can identify $\psi_1$ and $\psi_2$ with involutions on $E$, and since $\psi_1 \circ \psi_2 = \tau_1$ it follows that one of them is given by $T\mapsto T + R$ for some $R$ of order two.
Since $\psi_1$ and $\psi_2$ are different from $\tau_1\circ \tau_2$ it follows that $R\ne T_1+T_2$, thus since $T$ and $T+R$ correspond to the same curve, the claim follows.

\begin{figure}[htp]
	\centering
	\includegraphics[scale = 0.58]{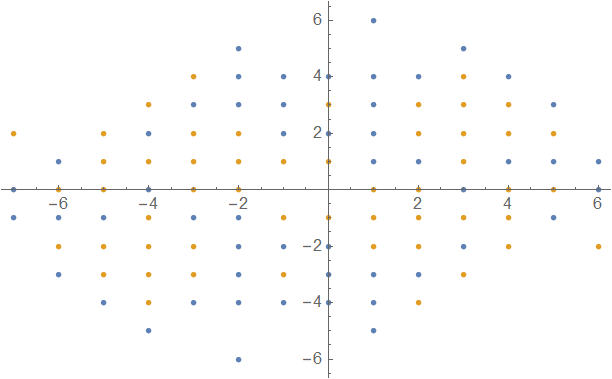}
	\caption{The blue (orange) point with coordinates $(n,m)$ represents the elliptic curve with root number one (minus one) that corresponds to the point $n Q_1 + m Q_2$ from the first example.}
	\label{fig:c}
\end{figure}

In our second example, see Section \ref{sec:other examples}, the elliptic curves with torsion group  $\Z/ 8 \Z$ and rank at least $3$ over $\Q$  are parametrized by rational points on elliptic curve $Y^2 = X^3 - 105987 X + 11743634$ of rank $2$ and with rational $2$-torsion subgroup. Denote the generators of the free part of the Mordell-Weil group by $Q_1= (-77,-4410)$ and  $Q_2=(805,21168)$. The root numbers of elliptic curves corresponding to the points $n Q_1+m Q_2$ for small $n,m \in \Z$ are presented in Figure \ref{fig:c} (as in the previous example, the points which differ by the point of order two correspond to the isomorphic elliptic curves, so they are not included in the figure). There are $52$ curves with the root number $1$, and $52$ curves with the root number $-1$, which suggests that the root numbers are evenly distributed also in this family.  Note that in this case, the figure has a symmetry $(n,m) \leftrightarrow (-n-1,-m)$.

Two parameterizations of rank three elliptic curves with torsion group $\Z/ 2 \Z \times \Z/ 6 \Z $  by rank two elliptic curves from Section \ref{sec:z26} are contained in the intersection of the rank two families from Sections \ref{sec:8.1} and \ref{sec:8.3}. More precisely, these families are described by curves $
D_1: r^2 s - 24 r s^2 + 168 rs - 336 r + 360 s^2 - 2700 s + 5040=0$ and $D_2:r^2 s^2 - \frac{15}{2} r^2 s + 14 r^2 - 12rs^2 + 84rs - 168r + 90s=0$ where $r$ and $s$ are parameters of the families from Sections \ref{sec:8.1} and \ref{sec:8.3}. Since these two equations have the same discriminant with respect to $r$ and $s$, we conclude that these (genus one) curves are not only birationally equivalent but also that the set of $s$ (and $r$) coordinates of the rational points on both of these curves agree. Hence, these two parameterizations are equal and can be described as follows.

\begin{figure}[htp]
	\centering
	\includegraphics[scale = 0.35]{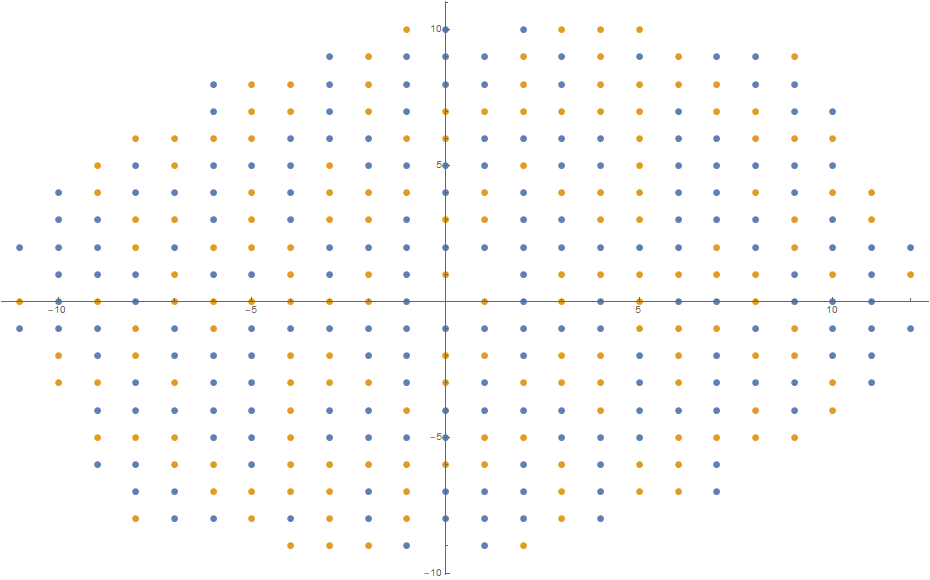}
	\caption{The blue (orange) point with coordinates $(n,m)$ represents the elliptic curve with root number one (minus one) that corresponds to the point $n R_1 + m R_2$.}
	\label{fig:d}
\end{figure}

The elliptic curves with torsion group  $\Z/ 2 \Z \times \Z/ 6 \Z $ and rank at least $3$ over $\Q$  are parameterized by rational points on elliptic curve $ y^2 = x^3 - x^2 - 456 x + 3456$ of rank $2$ and with rational $2$-torsion subgroup. Denote the generators of the free part of the Mordell-Weil group by $R_1= (20,-44)$ and  $R_2=(4/9, -1540/27)$. The root numbers of elliptic curves corresponding to the points $n R_1+m R_2$ for small $n,m \in \Z$ are presented in Figure \ref{fig:d} (as in the previous examples, the points which differ by the point of order two correspond to the isomorphic elliptic curves, so they are not included in the figure). There are $194$ curves with the root number $1$, and $168$ curves with the root number $-1$, which suggests that the root numbers are evenly distributed also in this family.  Note that in this case, the figure has a symmetry $(n,m) \leftrightarrow (-n+1,-m+1)$ which can be explained using involutions on curve  $D_1$ (or $D_2$) as in the first example.

\noindent{\it Final remark}.
Although our calculations of root numbers are limited to a relatively small number of cases,
it seems that they indicate that the heuristic in \cite{PPVW} needs some adjustments,
at least in the case of curves with torsion groups $\Z/ 8 \Z $ and $\Z/ 2 \Z \times \Z/ 6 \Z$.

\bigskip

{\bf Acknowledgement.} The authors would like to thank Ivica Gusi\'c and Maksym Voznyy for useful
comments on the previous version of this paper, and to the members of Mersenne Forum (https://mersenneforum.org/)  for the help with factorization in Section \ref{sec:heur}. The authors also would like to thank to the anonymous referees for several helpful suggestion.


 \end{document}